%% file: gen.tex
\magnification=\magstephalf
\def\to{\ \longrightarrow\ }

\def\nl{\hfill\break}

\input cprfonts

\input epsf
\overfullrule=0pt

\font\Bbb=msbm10

\font\secfont=cmbx10

\font\nam=cmr8
\font\aff=cmti8
\font\refe=cmr9

\mathchardef\square="0\hexa03
\def\qed{\hfill$\square$\par\rm}
\def\np{\vfill\eject}
\def\boxing#1{\ \lower 3.5pt\vbox{\vskip 3.5pt\hrule \hbox{\strut\vrule \ #1 \vrule} \hrule} }

\def\down#1{\ \lower 3.5pt\vbox{\vskip 3.5pt \hbox{\strut \ #1 \vrule} \hrule} }
\def\negdown#1{\ \lower 3.5pt\vbox{\vskip 3.5pt \hbox{\strut  \vrule \ #1 }\hrule} }
\def\adj{\hbox{\rm adj}}
\def\tr{\hbox{\rm tr}}

\hsize=6.3 truein
\vsize=9 truein

\baselineskip=13 pt
\parskip=\baselineskip
 1

\parindent=0pt

\def\R{\hbox{\Bbb R}}
\def\C{\hbox{\Bbb C}}

\def\d{\hbox{\rm d}}
\def\a{\hbox{\bf a}}
\def\b{\hbox{\bf b}}
\def\c{\hbox{\bf c}}
\def\x{\hbox{\bf x}}







\newif \iftitlepage \titlepagetrue

\def\diagram{\global\advance\diagramnumber by 1
$$\epsfbox{genfig.\number\diagramnumber}$$}
\def\ddiagram{\global\advance\diagramnumber by 1
\epsfbox{genfig.\number\diagramnumber}}

\newcount\diagramnumber
\diagramnumber=0

\newcount\secnum \secnum=0
\newcount\subsecnum
\newcount\defnum
\def\section#1{
                \vskip 10 pt
                \advance\secnum by 1 \subsecnum=0
                \leftline{\secfont \the\secnum \quad#1}
                }

\def\subsection#1{
                \vskip 10 pt
                \advance\subsecnum by 1 
                \defnum=1
                \leftline{\secfont \the\secnum.\the\subsecnum\ \quad #1}
                }

\def\definition{
                \advance\defnum by 1 
                \bf Definition 
\the\secnum .\the\defnum \rm \ 
                }

\def\lemma#1{
                \advance\defnum by 1 
                \par\bf Lemma  \the\secnum
.\the\defnum \sl \ #1 \par\rm
                }

\def\theorem#1{
                \advance\defnum by 1 
                \par\bf Theorem  \the\secnum
.\the\defnum \sl \ #1 \par\rm
               }
\def\cor#1{
                \advance\defnum by 1 
                \par\bf Corollary  \the\secnum
.\the\defnum \sl \ 
               }

\def\cite#1{
				\secfont [#1]
				\rm
}

\vglue 20 pt
\def\today{\ifcase\month\or
  January\or February\or March\or April\or May\or June\or
  July\or August\or September\or October\or November\or December\fi
  \space\number\day, \number\year}

\today

\centerline{\secfont GENERALIZED QUATERNIONS and}
\centerline{\secfont INVARIANTS of VIRTUAL KNOTS and LINKS}

\bigskip

\bigskip

\centerline{\nam ROGER FENN}
\centerline{\aff School of Mathematical Sciences, University of Sussex}
\centerline{\aff Falmer, Brighton, BN1 9RH, England}
\centerline{\aff e-mail addresses: rogerf@sussex.ac.uk}

\bigskip

\baselineskip=10 pt
\parskip=0 pt
\centerline{\nam ABSTRACT}
\leftskip=0.25 in
\rightskip=0.25in

{\nam In this paper we show how generalized quaternions including
$2\times 2$ matrices can be used to find solutions of the equation
$$[B,(A-1)(A,B)]=0.$$ These solutions can then be used to find
polynomial invariants of virtual knots and links.}

\leftskip=0 in
\rightskip=0 in
\baselineskip=13 pt
\parskip=\baselineskip

\parskip=\baselineskip
\section{Introduction}

Consider the algebra with the following presentation
$${\cal F}=\{A,B\mid A^{-1}B^{-1}AB-BA^{-1}B^{-1}A=B^{-1}AB-A\}.$$ In
this paper we will call this the {\bf fundamental} algebra and the
single relation will be called the fundamental relation or
equation. This relation arises naturally from attempts to find
representations of the braid group. Representations of the fundamental
algebra as matrices can be used to define representations of the
virtual braid group and invariants of virtual knots and links. In
\cite{BuF} we found a complete set of conditions for two classic
quaternions, $A, B$ to be solutions of the fundamental equation. In
this paper this result is generalised to give necessary and sufficient
conditions for generalized quaternions to satisfy the fundamental
relation, except in the case of all $2\times2$ matrices where only
sufficient conditions are given. Particularly, we define two
4-variable polynomials of virtual knots and links. In addition, we
give conclusive proof of the fact, only hinted at in earlier papers,
that invariants defined in this manner do not give any new invariants
for classical knots and links.

We are grateful to Jose Montesinos for suggesting the use of
generalized quaternions and to Steve Budden, Daan Krammer, Dale
Rolfsen and Bruce Westfield for helpful comments.

\np
\section{The fundamental equation  and its justification}
Given a set $X$ let $S$ be an endomorphism of $X^2$. In \cite{FJK},
such an $S$ is
called a {\bf switch} if
{\parindent=20pt
\item{1} $S$ is invertible and 
\item{2} the set theoretic Yang-Baxter equation
$$ (S\times id)(id\times
S)(S\times id)=(id\times S)(S\times id)(id\times S)$$
is satisfied. }
Switches are used in \cite{FJK} to define biracks and biquandles by the formula
$$S(a,b)=(b_a,a^b).$$

Given a switch $S$ there is a representation of the braid group $B_n$
into the group of permutations of $X^n$ defined by
$$\sigma_i\to (id)^{i-1}\times S\times (id)^{n-i-1}$$
where $\sigma_i$ are the standard generators.  Denote this representation by
$\rho=\rho(S, n)$. 

In this paper we will only be interested
in linear switches. So let $R$ be an associative
but not necessarily commutative ring and let $X$ be
a left $R$-module. Suppose
$$S=\pmatrix{A & B \cr C & D\cr}$$
where the matrix entries $A, B, C, D$ are elements of
$R$.
The 3$\times$3 matrices of the Yang-Baxter equation are
$$S\times id=\pmatrix{A&B&0\cr C&D&0\cr 0&0&1\cr}\quad
id\times S=\pmatrix{1&0&0\cr 0&A&B\cr 0&C&D\cr}$$
The representation $\rho$ is now into $n\times n$ matrices
with entries from $R$. 

Let us consider methods to find such switches $S$.
It is not difficult to see that the following 7 equations are
necessary and sufficient conditions for an invertible $S$ to be
a switch,
$$\matrix{%
1: A=A^2+BAC\hfill  &2: [B,A]= BAD\hfill \cr 
3: [C,D]= CDA\hfill  &4: D=D^2+CDB\hfill  \cr 
5: [A,C]= DAC\hfill  &6: [D,B]= ADB \hfill \cr 
\hfill 7: [C,B]=& ADA-DAD \hfill \cr}$$
where $[X,Y]=XY-YX$.

Examples of switches are
$$\leqalignno{
\hbox{ The identity }&&0:\cr
&&\cr
S=\pmatrix{0 & B \cr C & 1-BC\cr}& \hbox{ or }
S=\pmatrix{1-BC & B \cr C & 0\cr}&1:\cr}$$
where $B$ and $C$ are arbitary commuting invertible elements.
This is called the {\bf Alexander} switch. A special case
of this, when $B=1$, is called the {\bf Burau} switch.
$$\leqalignno{
S=\pmatrix{A & B \cr C & D\cr}&&2:\cr}$$ where
$A, A-1, B$ are invertible, $A, B$ do not commute and satisfy the
fundamental equation
$$A^{-1}B^{-1}AB-BA^{-1}B^{-1}A=B^{-1}AB-A$$
moreover
$$C=A^{-1}B^{-1}A(1-A),\ D=1-A^{-1}B^{-1}AB.$$
We will call this the {\bf non-commuting} switch.  A special case
of this is the matrix with quaternion entries
$$S=\pmatrix{1+i & j \cr -j & 1+i\cr}$$ called the
{\bf Budapest} switch.

If $S=\pmatrix{A & B \cr C & D\cr}$ is a switch then so is
$S(t)=\pmatrix{A & tB \cr t^{-1}C & D\cr}$ where $t$ is a commuting
variable. We say that $S(t)$ is $S$ {\bf augmented} by $t$.

In \cite{BuF} and \cite{BF} the following results can be found.
\theorem{Suppose $R$ is a division ring. Then any switch is one of
the examples above.}
\qed

Of course other types are possible, see \cite{Cs} in which divisors
of zero are used. 

The representation of the braid group induced by any non-commuting
switch looks complicated but is in fact equivalent to the Burau
representation. This has been pointed out previously by Dehornoy, see
\cite{De}. The following lemma gives an explicit proof.

\lemma{Let $S=\pmatrix{A & B \cr C & D\cr}$ be a non-commuting switch and
let $S'=\pmatrix{0 & 1 \cr Q & 1-Q\cr}$ be the Burau switch where
$Q=(1-A)(1-D)$.
Let $M$ be the $n\times n$ matrix
$$M=\pmatrix{
1&0&0&\cdots&0\cr
A&B&0&\cdots&0\cr
A&BA&B^2&\cdots&0\cr
&\vdots&&\vdots&0\cr
A&BA&B^2A&\cdots&B^{n-1}\cr}
$$
In words: the rows of $M$, after the first, start with $A$ and then the
previous row multiplied on the left by $B$. Clearly $M$ is invertible.
Then $\rho(S, n)=M^{-1}\rho(S', n)M$.}

{\bf Proof}
A calculation shows that $M\rho(S, n)=\rho(S', n)M$. In this calculation the
fundamental relation is used. For example $Q$ commutes with $B$. So to
prove that
$$B^iA^2+B^{i+1}C=QB^i+(1-Q)B^iA$$
we need to show that
$$A^2+BC=Q+(1-Q)A$$
which follows from the fundamental relation.
\qed

However, if we extend the representation to the {\bf virtual braid group},
defined below, then we get a representation which is {\bf not} equivalent
to the Burau.

The virtual braid group, $VB_n$ \cite{KK}, has generators
$\sigma_i,\ i=1,\ldots,n-1$ and braid group relations
$$\leqalignno{
\sigma_i \sigma_j&= \sigma_j \sigma_i, \qquad |i-j|>1 &i)\cr
\sigma_i \sigma_{i+1} \sigma_i&=\sigma_{i+1} \sigma_i
\sigma_{i+1}\hfill \cr }$$
In addition there are generators $\tau_i,\ i=1,\ldots,n-1$
and permutation group relations
$$\leqalignno{
{\tau_i}^2&=1 \hfill &ii) \cr
\tau_i \tau_j&= \tau_j \tau_i ,\qquad |i-j|>1 \cr
\tau_i \tau_{i+1} \tau_i&=\tau_{i+1} \tau_i \tau_{i+1}\hfill \cr } $$
together with mixed relations
$$\leqalignno{ \sigma_i \tau_j&= \tau_j \sigma_i ,\qquad
|i-j|>1 &iii)\cr
\sigma_i \tau_{i+1} \tau_i&=\tau_{i+1} \tau_i \sigma_{i+1}\hfill
\cr } $$

We can extend the representation $\rho(S, n)$ by sending the generator
$\tau_i$ to $(id)^{i-1}\times T \times (id)^{n-i-1}$ where
$T=\pmatrix{0 & 1 \cr 1 & 0\cr}$, that is, Burau with unit variable.

Consider now the element
$\beta=\sigma_2\sigma_1\tau_2\sigma_1^{-1}\sigma_2^{-1}\tau_1$ in
$VB_3$. If $S=\pmatrix{0 & B \cr C & 1-BC\cr}$, the Alexander switch,
then
$$\rho(\beta)=\pmatrix{1&0&(C^{-1}-B)(B-1)\cr 0&1&(C^{-1}-B)(1-B)\cr
0&0&1\cr}$$ So if $B\ne 1$ this representation can not be equivalent
to the Burau representation which has $B=1$.

We can now ask the following question: {\sl Let $K\subset B_n$ be the
kernel of the Burau representation. Let $\overline{K}$ denote the
normal closure of $K$ in $VB_n$. If $\beta$ is a virtual braid and
$\rho(\beta)=1$ for all switches $S$, is it true that $\beta$ lies in
$\overline{K}$?}

\section{Quaternion Algebras}
It is clear from the previous section that it is important to find
solutions to the fundamental equation. The main result in the next
section is a sufficient condition for two generalised quaternions to
satisfy the fundamental equation. Except for $2\times 2$ matrices,
this condition is also necessary.

In this section we describe the necessary algebra.  The results which
are already in the literature are mainly presented without proof.  For
more details see \cite{L}.

Let $F$ be a field of characteristic not equal to 2. Pick two non-zero
elements $\lambda,\ \mu$ in $F$. Let
$\bigl({{\lambda,\ \mu}\over F}\bigr)$
denote the algebra of dimension 4 over $F$ with basis $\{1,i,j,k\}$
and relations $i^2=\lambda,\ j^2=\mu,\ ij=-ji=k$.  The multiplication
table is given by
$$\bordermatrix{&i&j&k\cr i&\lambda&k&\lambda j\cr j&-k&\mu&-\mu i\cr
k&-\lambda j&\mu i&-\lambda\mu\cr}.$$
Throughout the paper a general quaternion algebra will be denoted by
${\cal Q}$.  Elements of ${\cal Q}$ are called (generalized)
quaternions. The field $F$ is called the {\bf underlying} field and
the elements $\lambda\ \mu$ the {\bf parameters} of the algebra. We
will denote quaternions by capital roman letters such as $A, B,
\ldots$ and (if pure) by bold face lower case, $\a,\b,\ldots$.  Field
elements, (scalars) will be denoted by lower case roman letters such
as $a, b, \ldots$ and lower case greek letters such as $\alpha, \beta,
\ldots$.

The classical quaternions are $\bigl({{-1,\ -1}\over \R}\bigr)$.
The algebra of 
$2\times 2$ matrices with entries in $F$ is $M_2(F)=\bigl({{-1,\ 1}\over F}
\bigr)$.

\subsection{Conjugation, Norm and Trace}

Let $A=a_0+a_1i+a_2j+a_3k$ be a quaternion where $a_0, a_1, a_2, a_3
\in F$. The coordinate $a_0$ is called the {\bf scalar} part of $A$
and the 3-vector $\a=a_1i+a_2j+a_3k$ is called the {\bf pure} part of
$A$. Evidently $A=a_0+\a$ is the sum of its scalar and pure parts and
is pure if its scalar part is zero and is a scalar if its pure part is
zero.

The {\bf conjugate} of $A$ is $\overline{A}=a_0-\a$, the {\bf norm} of $A$ is $N(A)=A\overline{A}$ and the 
{\bf trace} of $A$ is $\tr(A)=A+\overline{A}$. 

Conjugation is an anti-isomorphism of order 2. That is it satisfies
$$\overline{A+B}=\overline{A}+\overline{B},\quad
\overline{AB}=\overline{B}\;\overline{A},
\quad\overline{aA}=a\overline{A},\quad \overline{\overline{A}}=A.$$
Also $\overline{A}=A$ if and only if $A$ is a scalar and
$\overline{A}=-A$ if and only if $A$ is pure.

The norm is a scalar satisfying $N(AB)=N(A)N(B)$. We will denote the
set of values of the norm function by ${\cal N}$. It is a
multiplicatively closed subset of $F$ and ${\cal N}^*={\cal N}-\{0\}$
is a multiplicative subgroup of $F^*$. An element $A$ has an inverse
if and only if $N(A)\ne 0$ in which case
$A^{-1}=N(A)^{-1}\overline{A}$.

The trace of a quaternion is twice its scalar part.
\np
\subsection{Multiplying Quaternions}

Let $A, B$ be two quaternions.There is a bilinear form given by 
$$A\cdot B={1\over 2}(A\overline{B}+B\overline{A})={1\over 2}
(\overline{A}B+\overline{B}A)={1\over 2}\tr(A\overline{B}).$$
In terms of coordinates this is
$$A\cdot B=a_0b_0-\lambda a_1b_1-\mu a_2b_2+\lambda\mu a_3b_3.$$
Since $\lambda$ and $\mu$ are non-zero this is a non-degenerate form.
The corresponding quadratic form is
$$N(A)=a_0^2-\lambda a_1^2-\mu a_2^2+\lambda\mu a_3^2.$$
Let $\a, \b$ be pure quaternions. Then
$$\a\b=-\a\cdot\b+\a\times\b$$
where $$\a\cdot\b=-\lambda a_1b_1-\mu a_2b_2+\lambda\mu a_3b_3$$ is
the restriction of the bilinear form to the pure quaternions and
$\a\times\b$ is the {\bf cross product} defined symbolically by
$$\a\times\b=\Biggl|
\matrix{-\mu i&-\lambda j&k\cr a_1&a_2&a_3\cr b_1&b_2&b_3\cr}\Biggr|$$
The cross product has the usual rules of bilinearity and skew symmetry. The triple cross product expansion
$$\a\times(\b\times\c)=(\c\cdot\a)\b-(\b\cdot\a)\c$$
is easily verified.
The {\bf scalar triple product} is 
$$[\a,\b,\c]=\a\cdot(\b\times\c)=\lambda\mu
\Biggl|
\matrix{a_1&a_2&a_3\cr b_1&b_2&b_3\cr c_1&c_2&c_3\cr}\Biggr|$$
from which all the usual rules (except volume) can be deduced.

\subsection{Dependancy Criteria}

In this subsection we will consider conditions for sets of quaternions
to be linearly dependant or otherwise.  A non-zero element, $A$, of
${\cal Q}$ is called {\bf isotropic} if $N(A)=0$ and {\bf anisotropic}
otherwise. So only non-zero anisotropic elements have inverses. We
note the following theorem.
\theorem{The following statements about a quaternion algebra
${\cal Q}$ are equivalent.
{\parindent=20pt
\item{1.} ${\cal Q}$ contains an isotropic element.
\item{2.} ${\cal Q}$ is the sum of two hyperbolic planes.
\item{3.} ${\cal Q}$ is not a division algebra.
\item{4.} ${\cal Q}$ is $M_2(F)$.}
}

{\bf Proof} See \cite{L} p 58.

We will call a quaternion algebra above {\bf hyperbolic}. Otherwise it
is called {\bf anisotropic}. The classic quaternions are anisotropic:
$2\times2$ matrices are hyperbolic.

\lemma{A pair of pure quaternions $\a, \b$ is linearly dependant if and
only if $\a\times \b=0$.}

{\bf Proof} The proof is clear one way using the antisymmetry of the
cross product. Conversely suppose $\a\times \b=0$.  Then $(\a\times
\b)\times\c=(\a\cdot\c)\b-(\b\cdot\c)\a=0$. This can be made into a
linear dependancy by a suitable choice of $\c$, for example if
$\a\cdot\c\ne0$.
\qed
As a corollary we have the following
\lemma{Two quaternions commute if and only their pure parts are linearly
dependant.} \qed
Now we look for conditions for the triple of pure quaternions,
$\a,\b,\a\times \b$, to be linearly dependant.  The required condition
is given by the following lemma.

\lemma{The pure quaternions $\a,\b,\a\times \b$, are linearly dependant if
and only if 
$$N(\a)N(\b)=(\a\cdot\b)^2.$$ This is equivalent to the equations
$$N(\a\times\b)=-\mu(a_2b_3-a_3b_2)^2-\lambda(a_1b_3-a_3b_1)^2+
(a_1b_2-a_2b_1)^2=0,$$
ie $\a\times\b$ is isotropic or zero.}
{\bf Proof}
Three 3-dimensional vectors are linearly dependant if and only if the
determinant they form by rows is zero. In the case of pure quaternions
this means the scalar triple product is zero
$$[\a,\b,\c]=\a\cdot(\b\times\c)=\lambda\mu
\Biggl|
\matrix{a_1&a_2&a_3\cr b_1&b_2&b_3\cr c_1&c_2&c_3\cr}\Biggr|=0.$$
Replacing $\c$ with $\a\times\b$ and expanding out using the triple
cross product formula gives the first equation. Using the expansion formul\ae\ 
$$N(\a\times\b)=N(\a)N(\b)-(\a\cdot\b)^2$$ 
gives the second formula.
\qed
We have the following corollary.

\lemma{If $\a, \b$ are linearly independant pure quaternions and $\a\times \b$
is anisotropic, then the triple $\a,\b,\a\times \b$, is linearly independant.}

\subsection{$2\times2$ matrices}
We will interpret all the previous results in terms of $2\times 2$
matrices, $M_2(F)=\bigl({{-1,\ 1}\over F}\bigr)$. This is the only quaternion
algebra with zero divisors.

The generators of $\bigl({{-1,\ 1}\over F}\bigr)$ are, together with
the identity, the Pauli matrices
$$i=\pmatrix{0&1\cr -1&0\cr},\ j=\pmatrix{0&1\cr 1&0\cr},\ 
k=\pmatrix{1&0\cr 0&-1\cr}.$$
By an abuse of notation we will often confuse the scalar matrix
$\pmatrix{\nu&0\cr 0&\nu\cr}$ with the corresponding field element
$\nu$.

A general matrix can be written uniquely as
$$\pmatrix{\alpha&\beta\cr \gamma&\delta\cr}={1\over2}\bigl[(\alpha+\delta)+
(\beta-\gamma)i
+(\beta+\gamma)j+(\alpha-\delta)k\bigr]$$
Conversely
$$A=a_0+a_1i+a_2j+a_3k=\pmatrix{a_0+a_3&a_2+a_1\cr a_2-a_1&a_0-a_3\cr}$$

Conjugation is
$$\overline{A}=\adj A=\pmatrix{\delta&-\beta\cr -\gamma&\alpha\cr}=
\pmatrix{a_0-a_3&-a_2-a_1\cr a_1-a_2&a_0+a_3\cr}$$
and norm is
$$N(A)=\det A=\alpha\delta-\beta\gamma=a_0^2+a_1^2-a_2^2-a_3^2$$

The scalar part of $A$ is $a_0=\tr A/2=(\alpha+\delta)/2$ and the pure
part is
$$\pmatrix{a_3&a_2+a_1\cr a_2-a_1&-a_3\cr}=\pmatrix{(\alpha-\delta)/2&\beta\cr 
\gamma&(\delta-\alpha)/2\cr}$$

\subsection{Multiplying Matrices}

\lemma{Suppose $A, B\in M_2(F)=\left({{-1,\ 1}\over F}\right)$. Then $AB=AB$.}

The statement is deliberately provocative. It says that multiplying
$A, B$ as matrices and as quaternions is the same. This can be checked
directly. \qed

The above lemma allows quick checking of formul\ae\ so if
$$A=\pmatrix{\alpha_1&\alpha_2\cr 
\alpha_3&\alpha_4\cr}\hbox{ and }B=\pmatrix{\beta_1&\beta_2\cr
\beta_3&\beta_4\cr}\hbox{ then }
A\cdot B={1\over2}(\alpha_1\beta_4-\alpha_2\beta_3-\alpha_3\beta_2+
\alpha_4\beta_1)$$
If $\a=\pmatrix{\alpha_1&\alpha_2\cr \alpha_3&-\alpha_1\cr}$ and
$\b=\pmatrix{\beta_1&\beta_2\cr
\beta_3&-\beta_1\cr}$ are pure then
$\a\cdot\b=-\alpha_1\beta_1-(\alpha_2\beta_3+\alpha_3\beta_2)/2$
and
$$\a\times\b=\pmatrix{(\alpha_2\beta_3-\alpha_3\beta_2)/2&\alpha_1\beta_2-
\alpha_2\beta_1
\cr\alpha_3\beta_1-\alpha_1\beta_3 &(\alpha_3\beta_2-\alpha_2\beta_3)/2\cr}$$

\section{Solving the Fundamental Equation}

In this section we give a sufficient condition for two generalised
quaternions to satisfy the fundamental equation. Except for $2\times
2$ matrices, this condition is also necessary.

Let $A=a_0+{\bf a}$ and $B=b_0+{\bf b}$ be quaternions. We will need
the following easily checked lemma.

\lemma{Conjugation by multiplication is $B^{-1}AB=N(B)^{-1}\overline{B}AB$
where
$$\overline{B}AB=a_0(b_0^2+N(\b))+(b_0^2-N(\b))\a+2(\a\cdot\b)\b
+2b_0(\a\times\b).$$
The two commutators are
$$[A,B]=AB-BA=2\a\times\b$$
and $(A,B)=A^{-1}B^{-1}AB=N(A)^{-1}N(B)^{-1}\overline{A}\;\overline{B}AB$ where
$$\eqalign{\overline{A}\;\overline{B}AB
&=a_0^2b_0^2+b_0^2N(\a)+a_0^2N(\b)-N(\a)N(\b)+2(\a\cdot\b)^2\cr
&\qquad-2(b_0(\a\cdot\b)+a_0N(\b))\a
+2(a_0(\a\cdot\b)+b_0N(\a))\b+2(a_0b_0-\a\cdot\b)\a\times\b\cr}.$$}
\qed

We will call a quaternion, $A$, {\bf balanced} if $\tr(A)=N(A)\ne0$. A
quaternion $A=a_0+a_1i+a_2j+a_3k$ in a quaternion algebra with
parameters $\lambda, \mu$ is balanced if it lies on the quadric 3-fold
$$(a_0-1)^2-\lambda a_1^2-\mu a_2^2+\lambda\mu a_3^2=1,\quad a_0\ne0.$$
A balanced classical quaternion $A$ lies on the 3-sphere centre 1 and
radius 1. Note that if $A$ is balanced then $N(A-1)=1$. A pair of
invertible non-commuting quaternions, $A,B$, will be called {\bf
matching} if $A$ is balanced and $A\cdot B=0$.
\theorem{If the quaternion algebra is anisotropic then a necessary and sufficient condition for
the non-commuting, invertible quaternions $A, B$ to be solutions of
the fundamental equation is that they are a matching pair. Otherwise
the condition is only sufficient.}

{\bf Proof } The proof formally follows \cite{BuF}. In terms of
quaternions the equation is
$$\overline{A}\;\overline{B}AB-B\overline{A}\;\overline{B}A
=N(A)\overline{B}AB-N(A)N(B)A.$$
Using the formul\ae\ and notation developed above, the left hand side is
$$-4a_0N(\b)\a+4a_0(\a\cdot\b)\b-4(\a\cdot\b)\a\times\b$$
whereas the right hand side is
$$-2(a_0^2+N(\a))N(\b)\a+2(a_0^2+N(\a))(\a\cdot\b)\b
+2b_0(a_0^2+N(\a))\a\times\b$$
Considering half the difference of the two sides we arrive at
$$\c=(\tr(A)-N(A))N(\b)\a+(N(A)-\tr(A))(\a\cdot\b)\b+
(b_0(N(A)-\tr(A))+2A\cdot B)\a\times\b$$
So the two sides are equal if $\c=0$.

If $A, B$ is a matching pair then $\c=0$. Conversely if $\a, \b,
\a\times\b$ are linearly independent then their coefficients will be
zero and this implies that $A, B$ is a matching pair. The bilinear
form is definite unless the algebra is $M_2(F)$ and so, except for
this case, $\a, \b, \a\times\b$ will be linearly independent. \qed

\subsection{The Non-definite Case}

The condition for linear dependancy can be satisfied for the case when
the quaternions
are $2\times2$ matrices, ie $\lambda=-1, \mu=1$.  Let us try
$$a_2b_3-a_3b_2=1,\ a_1b_3-a_3b_1=1,\ a_1b_2-a_2b_1=0$$
One solution is
$$\a=i+tj-k,\ \b=j\hbox{ and }\a\times\b=k-i$$
where $t$ is an arbitary field element. So
$$\a-t\b+\a\times\b=0$$

We will return to this case in a later paper.

\section{The general matching pair}
In this section we will use the results of the previous section to describe
the most general matching pair. That is $A,B$ are $2\times 2$ matrices with
entries in some field and satisfying

{\parindent=20pt
\item{1.} $\tr(A)=\det(A)$ 
\item{2.} $A\adj(B)+B\adj(A)=0$ }

Since $B$ can be multiplied by any non-zero scalar we may assume temporarily
that

{\parindent=20pt
\item{3.} $\det(B)=1$.}

We can conjugate the matrices $A,B$ to simplify matters. Consider the
following two cases:

{\bf Case 1}, $A$ has two distinct eigenvalues and is diagonal.
$$A=\pmatrix{ a&0\cr 0&a/(a-1)\cr}\hbox{ and }
B=\pmatrix{ b&c\cr (b^2+a-1)/c(1-a)&
b/(1-a)\cr}$$
where $a,b,c$ are general and $c, 1-a$ must be invertible (ie non
zero).

Inverses are given by
$$A^{-1}={\displaystyle 1\over\displaystyle  a}\pmatrix{
1&0\cr 0&a-1\cr}$$
so $a$ must also be invertible and
$$B^{-1}=\adj(B)=\pmatrix{
b/(1-a)&-c\cr -
(b^2+a-1)/c(1-a) &b\cr}.$$

The $4\times 4$ matrix $S$ is $\pmatrix{A&B\cr C&D\cr}$ where 
$$C=A^{-1}B^{-1}A-A^{-1}B^{-1}A^2=\pmatrix{
b&c/(1-a)^2\cr (1-a)(b^2+a-1)/c&b/(1-a)\cr}
$$ and
$$D=1-A^{-1}B^{-1}AB=\pmatrix{
(2-3a+ab^2+a^2-2b^2)/(1-a)^2&(a-2)bc/(1-a)^2\cr
(a-2)b(b^2+a-1)/c(1-a)&(2-3a+ab^2+a^2-2b^2)/(1-a)\cr}
$$
Call this switch {$\bf E_2$}

{\bf Case 2} $A$ has one eigenvalue and is lower triangular
$$A=\pmatrix{
2&0\cr x&2\cr}.$$
Then $B$ has the form,
$$B=\pmatrix{
y&z\cr (xyz -2y^2-2)/2z & 
(xz-2y)/2\cr}$$
where $x,y,z$ are general and
$2, z$ must be invertible.

Inverses are given by
$$A^{-1}={1\over 4}\pmatrix{
2&0\cr -x&2\cr}.$$
and
$$B^{-1}=\adj(B)=\pmatrix{
(xz-2y)/2&-z\cr
 (-xyz+2y^2+2)/2z&y\cr}$$

The $4\times 4$ matrix $S$ is $\pmatrix{A&B\cr C&D\cr}$ where 
$$C=A^{-1}B^{-1}A-A^{-1}B^{-1}A^2=
\pmatrix{
y+xz&z\cr -(x^2z^2+3xyz+2y^2+2)/2z&-xz/2-y\cr}$$ and
$$D=1-A^{-1}B^{-1}AB=\pmatrix{
xyz/2&xz^2/2\cr -x(2y^2+xyz-2)/4&-xz(xz+2y)/4\cr}$$
Call this switch {$\bf E_1$}
\subsection{The variable $t$}

Let us now multiply $B$ by the scalar $t$, ie the switch is augmented
by $t$. This means we have two possible switches: one dependant on
four variables $a, b, c, t$ ($E_2$) and one dependant on four variables $x, y,
z, t$ ($E_1$).

\section{Determinants over Quaternion Algebras}

In order to define workable invariants we consider in this section a
determinental function on the matrices in $M_n({\cal Q})$. That is
$n\times n$ matrices with entries in a quaternion algebra ${\cal Q}$.
For background reading see \cite{As}. In fact the invariants defined
later can also be defined for any solutions of the fundamental
equation over a ring with a determinant function satisfying the rules
listed below.

If $R$ is a commutative ring let $\det:M_n(R)\to R$ denote the usual
determinant.  The classic quaternions, ${\cal H}$, may be embedded as
a subalgebra of $M_2(\C)$ and determinants taken in the usual way.
Our aim is to generalize this.

Suppose ${\cal Q}$ has underlying field $F$ and parameters $\lambda,
\mu$. Let $\overline{F}$ denote the algebraic closure 
of $F$. Embed $i, j, k$ in $M_2(\overline{F})$ by
$$i=\pmatrix{0&\sqrt{-\lambda}\cr -\sqrt{-\lambda}&0\cr},\
j=\pmatrix{0&\sqrt{\mu}\cr 
\sqrt{\mu}&0\cr},\ k=\pmatrix{\sqrt{-\lambda\mu}&0\cr
0&-\sqrt{-\lambda\mu}\cr}.$$
Define $d:M_n({\cal Q})\to \overline{F}$ as the composition of the
embedding $M_n({\cal Q})\subset M_{2n}(\overline{F})$ with $\det$.

Alternatively the determinant function may be defined by induction on
the size of the matrices. The value $d(A)=N(A)$ starts the
induction. Consider a matrix in $M_n({\cal Q})$. This may be reduced
to diagonal form, by multiplying on the left and the right by
elementary matrices having unit determinant, (see below). Suppose this
matrix has diagonal elements $d_1,\ldots,d_n$. Define the determinant
as $d=N(d_1)\cdots N(d_n)$. So the determinant function takes values
in ${\cal N}$. For $M_2(F)$ this subset is the whole of $F$: for
classic quaternions it is the non-negative reals.

The determinant function satisfies the rules
{\parindent=20pt
\item{0.} $d(M)=0$ if and only if $M$ is singular, moreover $d(MN)=d(M)d(N)$.
It follows that $d(1)=1$.
\item{1.} $d$ is unaltered by a permutation of the rows (columns). 
\item{2.} If a row (column) is multiplied on the left (right) by a unit
then $d$ is multiplied by $d$ of that unit.
\item{3.} $d(M)$ is unaltered by adding a left multiple of a row to another
row or a right multiple of a column to another column.
\item{4.} $d\pmatrix{x&{\bf u}\cr{\bf 0}&M\cr} =N(x)d(M)$ where ${\bf u}$
is any row vector and ${\bf 0}$ is a zero column vector.
\item{4'.} $d\pmatrix{x&{\bf 0}\cr{\bf v}&M\cr} =N(x)d(M)$ where ${\bf v}$
is any column vector and ${\bf 0}$ is a zero row vector.
\item{5.} $d(M^*)=\d(M)$ where $M^*=\overline{M^T}$ denotes the Hermitian conjugate.
\item{6.} if the entries in $M$ all commute then $d(M)=\det^2(M)$.}

An {\bf elementary} matrix of type 1 is a permutation matrix. An elementary
matrix of type 2 is the identity matrix with one diagonal entry replaced by a
unit and an elementary matrix of type 3 is a square matrix with zero entries
except for $1$'s down the diagonal and one other entry off diagonal.

The properties $i.$ above for $i=1, 2, 3$ follow from multiplying $M$ on
the right or left by an elementary matrix of type $i$.

The matrix $S$ can be written as a product of elementary matrices
$$\pmatrix{%
A & B \cr C & D\cr}=
\pmatrix{%
A & 0 \cr 0 & 1}
\pmatrix{%
1 & 0 \cr C & 1\cr}
\pmatrix{%
1 & 0 \cr 0 & 1-A^{-1}\cr}
\pmatrix{%
1 & A^{-1}B \cr 0 & 1}.
$$
Note that $1-A^{-1}$ is invertible.

Hence $d(S)=d(A)d(1-A^{-1})=d(A-1)=1$. Therefore the representation of $VB_n$,
induced by such an $S$, is into $SL(F,2n)$.

\section{Virtual Knots and Links}
Recall that classical knot theory can be described in terms of knot and link
diagrams. A {\bf diagram} is a 4-regular plane graph (with extra structure at
its nodes representing the crossings in the link) on a plane and
implicitly on a two-dimensional sphere $S^{2}$. 
We say that two such diagrams are {\bf equivalent} if there is a sequence of
moves of the types indicated in part (A) of Figure 1 (The Reidemeister Moves)
taking one diagram to the other. These moves are performed locally on the
4-regular plane graph (with extra structure) that constitutes the link
diagram.

Virtual knot theory is an extension of classical knot theory, see
\cite{K}. In this extension one adds a {\bf virtual crossing} (See
Figure 1) that is neither an over-crossing nor an under-crossing. We
shall refer to the usual diagrammatic crossings, that is those without
circles, as {\bf real} crossings to distinguish them from the virtual
crossings. A virtual crossing is represented by two crossing arcs with
a small circle placed around the crossing point. The arcs of the graph
joining adjacent classical crossings are called the {\bf semi-arcs} of
the diagram.

In addition to their application as a geometric realization of the
combinatorics of a Gauss code, virtual links have physical,
topological and homological applications. In particular, virtual links
may be taken to represent a particle in space and time which
dissappears and reappears. A virtual link may be represented, up to
stabilisation, by a link diagram on an orientable surface,
\cite{Ku}. If the surface has minimal genus then this representation
is unique. Finally an element of the second homology of a rack space
can be represented by a labelled virtual link, see \cite{FRS1}
\cite{FRS2}. Since the rack spaces form classifying spaces for
classical links the study of virtual links may give information about
classical knots and links.

The allowed moves on virtual diagrams are a generalization of the Reidemeister
moves for classical knot and link diagrams. We show the classical Reidemeister
moves as part (A) of Figure 1. These classical moves are part of virtual
equivalence where no changes are made to the virtual crossings. Taken by
themselves, the virtual crossings behave as diagrammatic permutations.
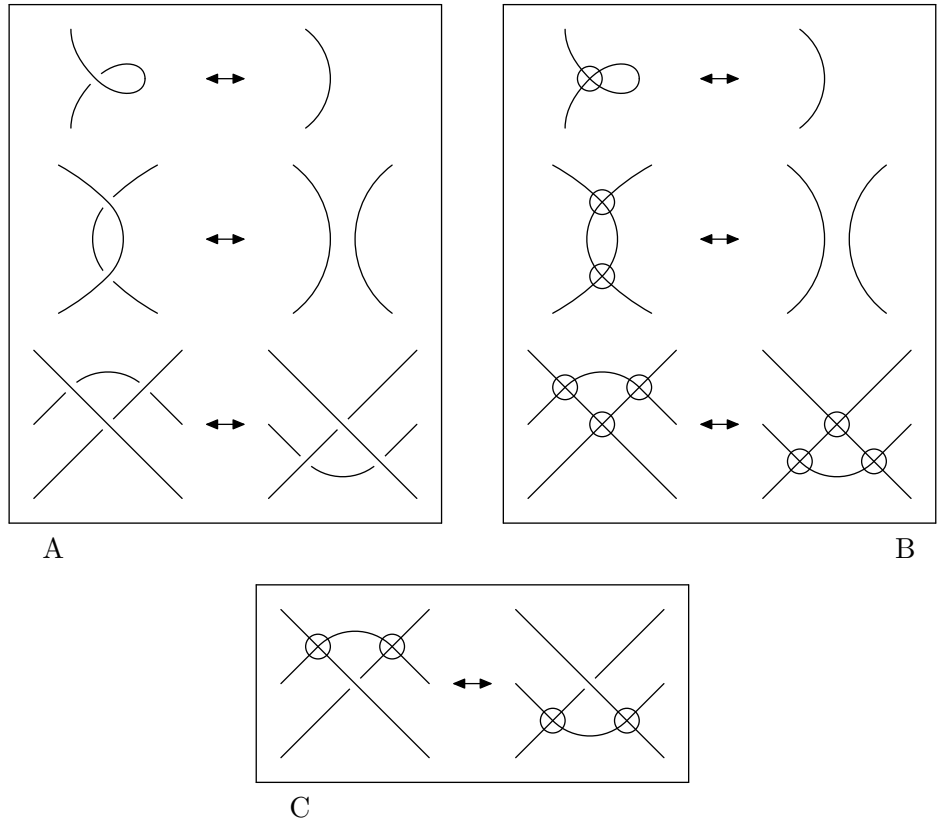
Specifically, we have the flat Reidemeister moves (B) for virtual crossings as
shown in Figure 1. In Figure 1 we also illustrate a basic move (C) that
interrelates real and virtual crossings. In this move an arc going through a
consecutive sequence of two virtual crossings can be moved across a single
real crossing. In fact, it is consequence of moves (B) and (C) for virtual
crossings that an arc going through any consecutive sequence of virtual
crossings can be moved anywhere in the diagram keeping the endpoints fixed and
writing the places where the moved arc now crosses the diagram as new virtual
crossings. This is shown schematically in Figure 2. We call the move in Figure
2 the {\bf detour}, and note that the detour move is equivalent to having all
the moves of type (B) and (C) of Figure 1. This extended set of moves
(Reidemeister moves plus the detour move or the equivalent moves (B) and (C))
constitutes the set of moves for diagrams of virtual knots and links.
\diagram \diagram
The topological interpretation for this virtual theory in terms of
embeddings of links in thickened surfaces is a useful idea. See
\cite{KK}, \cite{Ku}. Regard each virtual crossing as a shorthand for
a detour of one of the arcs in the crossing through a 1-handle that
has been attached to the 2-sphere of the original diagram. The two
choices for the 1-handle detour are homeomorphic to each other (as
abstract surfaces with boundary a circle) since there is no a priori
difference between the meridian and the longitude of a torus. By
interpreting each virtual crossing in this way, we obtain an embedding
of a collection of circles into a thickened surface $S_{g} \times \R$
where $g$ is the number of virtual crossings in the original diagram
$L$, $S_{g}$ is a compact oriented surface of genus $g$ and $\R$
denotes the real line. Thus to each virtual diagram $L$ we obtain an
embedded disjoint union of circles in $S_{g(L)} \times \R$ where
$g(L)$ is the number of virtual crossings of $L$.  We say that two
such surface embeddings are {\it stably equivalent} if one can be
obtained from another by isotopy in the thickened surfaces,
homeomorphisms of the surfaces and the addition or subtraction of
empty handles. Then we have the \smallbreak

\theorem{\sl Two virtual link diagrams
are equivalent if and only if their correspondent surface embeddings are
stably equivalent, \cite{KK}, \cite{Ku}. }

\smallbreak

The surface embedding interpretation of virtuals is useful since it
converts their equivalence to a topological question. The diagrammatic
version of virtuals embodies the stabilization in the detour moves. We
shall rely on the diagrammatic approach here.

\section{The invariant knot modules}

In this section, we shall begin to show how the previous algebra can
give rise to virtual knot invariants. Given an associative ring $R$, a
$2\times 2$ matrix $S$ with entries in $R$ and a virtual link diagram
$\cal D$, we define a presentation of an $R$-module which depends only
on the link class of the diagram provided $S$ is a switch. This
construction also works for classical knots and links but is only the
Alexander module in disguise. The generators are the semi-arcs of
$\cal D$, that is the portion of the diagram bounded by two adjacent
classical crossings. There are 2 relations for each classical
crossing.

Suppose the diagram $\cal D$ has $n$ classical crossings. Then there
are $2n$ semi-arcs labelled $a, b, \ldots$. These will be the
generators of the module.  Let the edges of a positive real crossing
in a diagram be arranged diagonally and called geographically {\it NW,
SW, NE} and {\it SE}. Assume that initially the crossing is oriented
and the edges oriented towards the crossing from left to right ie west
to east. The {\bf input} edges, oriented towards the crossing, are in
the west and the edges oriented away from the crossing, the {\bf
output} edges, are in the east. Let $a$ and $b$ be the labels of the
input edges with $a$ labelling SW and $b$ labelling NW. For a positive
crossing, $a$ will be the label of the undercrossing input and $b$ the
label of the overcrossing input. Suppose now that
$$S(a,b)^T=(c,d)^T\hbox{ where }S=\pmatrix{A&B\cr C&D\cr}.$$
Then we label the undercrossing output NE by $d$ and we label the
overcrossing output SE by $c$.

For a negative crossing the direction of labelling is reversed. So $a$
labels SE, $b$ labels NE, $c$ labels SW and $d$ labels NW.

Finally for a virtual crossing the labellings carry across the
strings. This corresponds to the twist function $T(a,b)=(b,a)$.

The following figure shows the labelling for the three kind of
crossings and the corresponding relations for the 2 classical crossings.
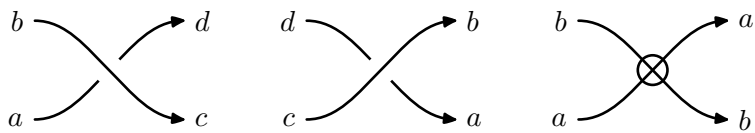
\diagram
\centerline{$c=Aa+Bb\quad d=Ca+Db$}

The diagram therefore gives rise to a presentation of an $R$-module
with $2n$ generators and $2n$ relations. Note that in all cases $B, C$
are invertible since the identity switch is uninteresting.
\theorem{The module defined above for any diagram $D$ is invariant
under the Reidemeister moves, and hence is a knot invariant, if $S$ is
a switch.}

{\bf Proof }
The proof of the above can be found in the papers \cite{FJK, BF, BuF},
For the
convenience of the reader we show how the module is invariant under
the Reidemeister moves.

Refering back to the picture of the two relations defined by a
crossing, it is convenient to think of the action from left to right
on a positive crossing as being the action of $S$ and the action from
right to left as being $S^{-1}$.

Consider the action from top to bottom as being $S^-_+$ and the action
from bottom to top as being $S^+_-$. By solving the equations of the
labellings we see that these matrices are
$$S^+_-=\pmatrix{%
DB^{-1}&C-DB^{-1}A\cr B^{-1}&-B^{-1}A\cr}\quad
S^-_+=\pmatrix{%
-C^{-1}D&C^{-1}\cr B-AC^{-1}D&AC^{-1}\cr}
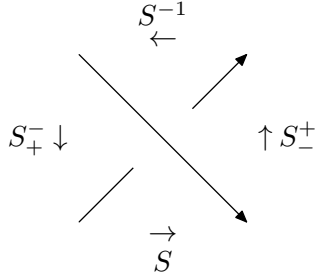
$$
\diagram
We call $S^+_-$ and $S^-_+$ the {\bf sideways} matrices. They are
invertible since $S$ is. Also $(S^{-1})^+_-=(S^-_+)^{-1}$ and
$(S^{-1})^-_+=(S^+_-)^{-1}$ and
$$S^+_-(a,a)=(\lambda a,\lambda a)\hbox{ and }
S^-_+(a,a)=(\lambda^{-1} a,\lambda^{-1} a)$$
where $\lambda=B^{-1}(1-A)=(1-D)^{-1}C$. So the sideways matrices
preserve the diagonal. This has the curious consequence that a linear
switch which is a birack is also a biquandle in the sense of
\cite{FJK}.

For a negative crossing the actions are equal
but with opposite orientation.

Assume for simplicity that we are dealing with a knot. The link case is
similar and details can safely be left to the reader.
We have a right $R$-module with a finite (square) presentation.

Following the orientation of
the knot, label the semi-arcs with $R$-variables $x_1, x_2, \ldots, x_{2n}$.
By an $R$-variable we mean a symbol standing in for any element of $R$.

At each crossing there is a relation of the form
$$\pmatrix{A&B\cr C&D\cr}\pmatrix{x_i\cr x_j\cr}=
\pmatrix{x_{j+1}\cr x_{i+1}\cr}
\hbox{ or }
\pmatrix{A&B\cr C&D\cr}\pmatrix{x_i\cr x_j\cr}=
\pmatrix{x_{j-1}\cr x_{i-1}\cr}$$
depending on whether the crossing is positive or negative. As is the usual
custom, indices are taken modulo $2n$.

The relations can now be written in matrix form as $M{\bf x}={\bf 0}$ where
$M$ is a $2n\times 2n$ matrix and ${\bf x}=(x_1,x_2 \ldots, x_{2n})^T$. The
non-zero entries in each row of the matrix are $A, B, -1$ or $C, D, -1$.

Let ${\cal M}={\cal M}(S,D)$ be the module defined by these relations. We now
show that the modules defined by diagrams representing the same virtual link
are isomorphic. We do this by showing that a single Reidemeister move defines
an isomorphism. The proof has the same structure as the proof, say, that the
Alexander module of a classical link is an invariant as in \cite{Alex} but we
give the details because of the care needed due to non-commutativity.

Any module defined by a presentation of the form $M{\bf x}={\bf 0}$ is
invariant under the following moves and their inverses applied to the matrix
$M$. 
 
{\parindent=20pt
\item{1.} permutations of rows and columns,
\item{2.} multiplying any row on the left or any column on the right by a
unit,
\item{3.} adding a left
multiple of a row to another row or a right multiple of a column to another
column, 
\item{4.} changing $M$ to $\pmatrix{x&{\bf u}\cr{\bf 0}&M\cr}$ where
$x$ is a unit, ${\bf u}$ is any row vector and ${\bf 0}$ is a zero column
vector,
\item{5.} repeating a row.
}

The operations $i.$ above for $i=1, 2, 3$ are equivalent to multiplying $M$ on
the right or left by an elementary matrix of type $i$.

Recall that the matrix $S$ can be written as a product of elementary matrices
$$\pmatrix{%
A & B \cr C & D\cr}=
\pmatrix{%
A & 0 \cr 0 & 1}
\pmatrix{%
1 & 0 \cr C & 1\cr}
\pmatrix{%
1 & 0 \cr 0 & 1-A^{-1}\cr}
\pmatrix{%
1 & A^{-1}B \cr 0 & 1}.
$$

Now consider the module ${\cal M}$ defined above. Clearly the presentation
is unaltered by any of the basic moves which involve the virtual crossing. So
we look to see the changes induced by the classical Reidemeister moves and
check that the presentation matrix $M$ is only changed by the above 5 moves.
Assume
$$M=\pmatrix{
m_{11}&\ldots&m_{1n-1}&m_{1n}\cr
m_{21}&\ldots&m_{2n-1}&m_{2n}\cr
\vdots&\ddots&\vdots&\vdots&\cr
m_{n1}&\ldots&m_{nn-1}&m_{nn}\cr
}$$

Firstly, consider a Reidemeister move of the first kind. 
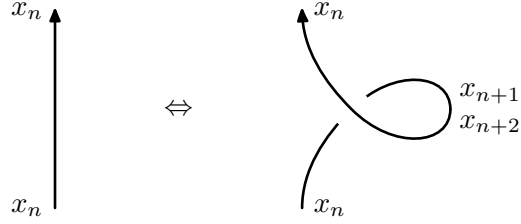
\diagram This
introduces (or deletes) two new equal generators $x_{n+1}=x_{n+2}$. Because
$S^-_+$ and $S^+_-$ preserve the diagonal, (the biquandle condition, see
\cite{FJK}) the output ($x_{n+3}$) is the same as the input ($x_{n}$). The
generator $x_{n+1}$ is equal to $\lambda^{-1}x_{n}$ where
$\lambda=B^{-1}(1-A)$.

So up to reordering of the columns the relation matrix is changed by
$$M\ \Leftrightarrow\ \pmatrix{
m_{11}&\ldots&m_{1n-1}&m_{1n}&0&0\cr
m_{21}&\ldots&m_{2n-1}&m_{2n}&0&0\cr
\vdots&\ddots&\vdots&\vdots&\vdots&\vdots\cr
m_{n1}&\ldots&m_{nn-1}&m_{nn}&0&0\cr
0&\ldots&0&0&1&-1\cr
0&\ldots&0&1&-\lambda&0\cr}
\Leftrightarrow\ \pmatrix{
m_{11}&\ldots&m_{1n-1}&m_{1n}&0&0\cr
m_{21}&\ldots&m_{2n-1}&m_{2n}&0&0\cr
\vdots&\ddots&\vdots&\vdots&\vdots&\vdots\cr
m_{n1}&\ldots&m_{nn-1}&m_{nn}&0&0\cr
0&\ldots&0&0&0&-1\cr
0&\ldots&0&1&\lambda&0\cr
}$$
Since $\lambda$ is a unit this does not alter the module.

There are other possible inversions and mirror images of the above which
can be dealt with in a similar fashion.

Secondly, consider a Reidemeister move of the second kind. 
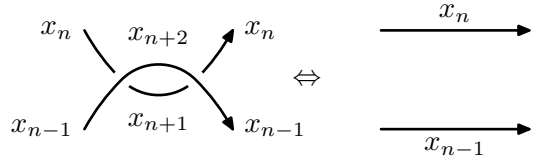
\diagram Again the
outputs are unchanged from the inputs $x_{n-1}, x_n$ because of the relation
$S^{-1}S=1$.

Two new generators $x_{n+1}$ and $x_{n+2}$ are introduced (or deleted).
They are related by the equations
$$x_{n-1}=Ax_{n+1}+Bx_{n+2}\hbox{ and }x_{n}=Cx_{n+1}+Dx_{n+2}.$$

This has the following effect on the relation matrix.
$$M\ \Leftrightarrow\ \pmatrix{
m_{11}&\ldots&m_{1n-1}&m_{1n}&0&0\cr
m_{21}&\ldots&m_{2n-1}&m_{2n}&0&0\cr
\vdots&\ddots&\vdots&\vdots&\vdots&\vdots\cr
m_{n1}&\ldots&m_{nn-1}&m_{nn}&0&0\cr
0&\ldots&0&-1&A&B\cr
0&\ldots&-1&0&C&D\cr
}.
$$
Since $S$ is a product of elementary matrices this does not alter the module.

The other possible inversions and mirror images of the above can be dealt with
in a similar fashion but it is worth looking at the case where the two arcs
run in opposite directions. The right outputs are unchanged from the left
inputs by the relation $S^+_-(S^{-1})^+_-=1$.
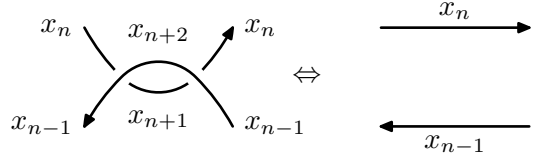
\diagram
The changes to the
relation matrix are given by $$M\ \Leftrightarrow\  \pmatrix{
m_{11}&\ldots&m_{1n-1}&m_{1n}&0&0\cr
m_{21}&\ldots&m_{2n-1}&m_{2n}&0&0\cr
\vdots&\ddots&\vdots&\vdots&\vdots&\vdots\cr
m_{n1}&\ldots&m_{nn-1}&m_{nn}&0&0\cr
0&\ldots&-1&A&B&0\cr
0&\ldots&0&C&D&-1\cr
}\ \Leftrightarrow\ 
\pmatrix{
m_{11}&\ldots&m_{1n-1}&m_{1n}&0&0\cr
m_{21}&\ldots&m_{2n-1}&m_{2n}&0&0\cr
\vdots&\ddots&\vdots&\vdots&\vdots&\vdots\cr
m_{n1}&\ldots&m_{nn-1}&m_{nn}&0&0\cr
0&\ldots&0&0&1&0\cr
0&\ldots&0&0&0&-1\cr
}
$$ 

Using the fact that $B$ is a  unit. This doesn't change the module.

Finally, consider a Reidemeister move of the third kind. 
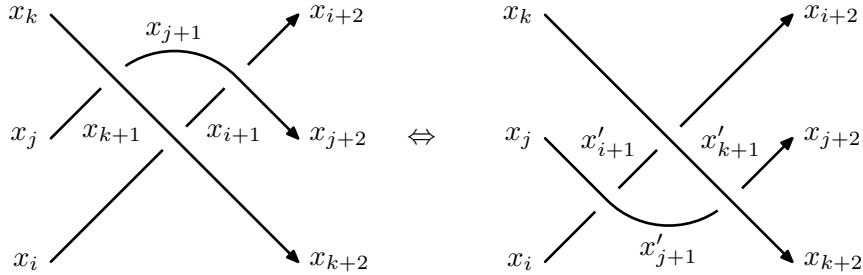
\diagram The outputs
$x_{i+2}, x_{j+2}, x_{k+2}$ are unaltered by the Reidemeister move because of
the Yang-Baxter equations. The inner generators $x_{i+1}, x_{j+1}, x_{k+1}$
are related to the inputs $x_{i}, x_{j}, x_{k}$ by the following matrix
$$\pmatrix{C&DA&DB\cr 0&C&D\cr 0&A&B\cr}$$ and the inner
generators $x'_{i+1}, x'_{j+1}, x'_{k+1}$ are related to $x_{i}, x_{j}, x_{k}$
by the following matrix $$\pmatrix{C&D&0\cr A&B&0\cr
AC&AD&B\cr}.$$ Both are the product of elementary matrices and the
proof is finished. \qed

\section{Determinant Invariants}

\subsection{The Determinant $\Delta_0$}
Given a module with a square presentation the obvious invariant of the
module is the determinant, if it can be defined. This will be the case
if the ring is represented by matrices with commuting entries, for
example the ring of generalised quaternions. In this case if $d$
denotes the determinant and $M\x =0$ is the presentation let
$\Delta_0=d(M)$.  Since the module depends on the switch $S$ we
illustrate this dependency by $\Delta_0=\Delta_0(S)$.

A close look at how the presentation of the module changes under the
Reidemeister moves shows that $\Delta_0$ is invariant up to
multiplication by $d(B)$ or $d(C)$. Typically $d(B)$ is denoted by the
variable $t$ and $d(C)$ is $t^{-1}$. If we take the switch to be $E_1$ ($E_2$)
defined in section 5  then $\Delta_0$ is a polynomial $p_1$ ($p_2$) in the four
variables $x, y, z, t$ ($a, b, c, t$). We can normalise these polynomial so
that as a polynomial in $t$ it has a non-zero constant term and only positive
powers of $t$.

Let us illustrate the previous discussion by calculating invariants for the
{\bf virtual trefoil} as shown in the figure.
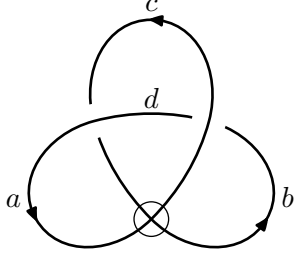
\diagram
If we label as
indicated then the module has a presentation with 4 generators
$a,b, c, d$ and relations $c=Ab+Ba,\ a=Ac+Bd,\ b=Cc+Dd,\ d=Cb+Da$.
Restricting to the $E_1$ case gives the polynomial $p_1$ equal to
$${64+4x^3tz^3+4x^3t^3z^3+128t^2-64xtz-x^4t^2z^4-64xt^3z+8x^2t^2z^2
-4x^2z^2+64t^4-4x^2t^4z^2\over 16}$$
For the $E_2$ case we get the polynomial $p_2$ equal to
$$
-{{\matrix{
-4+16\,a-40\,tb{a}^{3}-4\,{t}^{4}{a}^{4}+16\,{t}^{4}{a}^{3}-
24\,{t}^{4}{a}^{2}+16\,{t}^{4}a-36\,t{b}^{3}{a}^{2}+72\,{a}^{2}tb+4\,{
b}^{2}+13\,{a}^{2}{t}^{4}{b}^{2}\cr
+40\,t{b}^{3}a+4\,{t}^{4}{b}^{2}-4\,{t}^{4}-24\,{a}^{2}-12\,{b}^{2}a
+13\,{a}^{2}{b}^{2} +16\,{a}^{3}
-12\,{t}^{4}{b}^{2}a+16\,tb-6\,{t}^{4}{b}^{2}{a}^{3}\cr
+{t}^{4}{b}^{2}{a}^{4}-56\,
atb-16\,t{b}^{3}+{b}^{2}{a}^{4}-6\,{b}^{2}{a}^{3}+8\,b{a}^{4}t-2\,{b}^
{3}{a}^{4}t+14\,{b}^{3}{a}^{3}t-4\,{a}^{4}-8\,{t}^{2}{a}^{4}+32\,{t}^{
2}{a}^{3}\cr
-48\,{t}^{2}{a}^{2}+32\,{t}^{2}a-8\,{t}^{2}+
{a}^{4}{t}^{2}{b}
^{4}-2\,{a}^{4}{t}^{2}{b}^{2}-8\,{a}^{3}{t}^{2}{b}^{4}-8\,{t}^{2}{b}^{
2}+16\,{t}^{2}{b}^{4}+16\,{t}^{3}b-16\,{t}^{3}{b}^{3}-\cr
26\,{a}^{2}{t}^{
2}{b}^{2}+12\,{a}^{3}{t}^{2}{b}^{2}+24\,{a}^{2}{t}^{2}{b}^{4}+
8\,{a}^{
4}{t}^{3}b-40\,{a}^{3}{t}^{3}b+72\,{a}^{2}{t}^{3}b+24\,a{t}^{2}{b}^{2}
-32\,a{t}^{2}{b}^{4}\cr
-56\,a{t}^{3}b-2\,{t}^{3}{b}^{3}{a}^{4}+
14\,{t}^{3
}{b}^{3}{a}^{3}-36\,{t}^{3}{b}^{3}{a}^{2}+40\,{t}^{3}{b}^{3}a}\over
{t}^{2}\left( a-1 \right) ^{4}}}
$$

Note that the fundamental quandle (and hence group) as defined by
the Wirtinger presentation is trivial.

The following virtual knot is interesting in having a trivial
Jones-polynomial as well as a trivial fundamental rack.
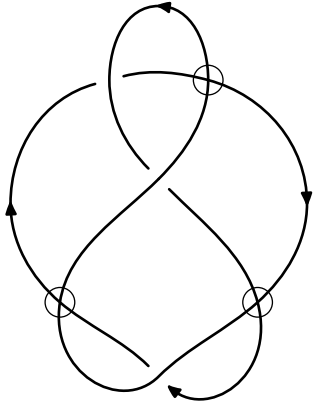
\diagram
In this case if $S$ is the Alexander switch then
$$\Delta_0=(B-1)(C^2(B+1)-C(B+1)(B^{-1}+1)-B).$$

\subsection{The Determinant $\Delta_1$}
For many knots and links, including
the classical, the determinant $\Delta_0$ is zero.

For example, as we have seen earlier any switch $S$ with entries in
the ring $R$ defines a representation of the virtual braid group
$VB_n$ into the group of invertible $n\times n$ matrices with entries
in $R$ by sending the standard generator $\sigma_i$ to
$S_i=(id)^{i-1}\times S
\times (id)^{n-i-1}$ and the generator $\tau_i$ to
$T_i=(id)^{i-1}\times T \times (id)^{n-i-1}$. This representation is
denoted by $\rho=\rho(S, n)$. For classical braids this representation
is equivalent to the Burau representation and so we would expect the
closure of a classical braid to have $\Delta_0$ zero. We now confirm
this by looking at the fixed points of $S_i$ both on the left and
right.

\lemma{Let $P=A^{-1}B^{-1}A$ and $Q=B^{-1}(1-A)$. Then
$$(P^{n-1},\ldots,P,1)S_i=(P^{n-1},\ldots,P,1)\qquad (*)$$
{\sl and}
$$S_i(1,Q,\ldots,Q^{n-1})^T=(1,Q,\ldots,Q^{n-1})^T\qquad (**).$$}

{\bf Proof } We need only check that
$$(P,1)S=(P,1)\hbox{ and } S(1,Q)^T=(1,Q)^T.$$
\qed

Therefore the following lemma gives a necessary condition for
the knot or link to be classical.

\theorem{For all classical knots and links $\Delta_0=0$.}

{\bf Proof } Since $\Delta_0$ is an invariant of the module
we can assume that the diagram from which it is defined is the closure
of a braid. However from $*$ (or $**$) there is a linear relationship
amongst the rows (columns), so $\Delta_0=0$.
\qed

The Kishino knots $K_1, K_2$ and $K_3$ are illustrated below.
$$\ddiagram\qquad\ddiagram\qquad\ddiagram$$
All are ways of forming the connected sum of two unknots. $K_1$ and
$K_2$ are mirror images and $K_3$ is amphich\ae ral. Both have trivial
racks and Jones polynomial. The invariant $\Delta_0$ is zero in all
three cases.

It is clear that we need an invariant for these cases. Let $M$ be the
$n\times n$ presentation matrix. Let $M_1, M_2, \ldots, M_{n^2}$ be
the submatrices obtained by deleting a row and a column. Let $d_1,
d_2, \ldots, d_{n^2}$ be the determinants. These all lie in a ring of
polynomials with coefficients in a field. Therefore the determinants
have an hcf, call it $\Delta_1$, which is well defined up to
multiplication by a unit. Now look at what happens to this
construction under a Reidemeister move. The hcf, $\Delta_1$, is
multiplied by $d$ of a unit.

Returning to the Kishino knots, a calculation with the Alexander switch
shows that for $K_1$, $\Delta_1$ is $1+B-CB$ and for $K_2$,
$\Delta_1$ is $1+C-CB$. Since these are neither units nor associates
in the ring, $K_1, K_2$ are non-trivial and non amphich\ae ral.

On the other hand for $K_3$, $\Delta_1$ is 1. We will show shortly 
that $K_3$
is non-trivial by using the Budapest switch augmented by $t$. Then
$\Delta_1=2+5t^2+2t^4$ see \cite{BuF}.

For a classical knot or link the invariant $\Delta_1$ is not just the
Alexander polynomial in disguised form but is independant of the
deleted rows or columns chosen, up to multiplication by a power of $t$.

\theorem{ Let ${\cal D}$ be a diagram of a classical knot or link. If
$M$ is the presentation matrix associated with ${\cal D}$.  Let
$\Delta_1=d(M_{ij})$ where $M_{ij}$ is obtained from $M$ by deleting
the $i$th row and the $j$th column. Then $\Delta_1$ is independant of
$i, j$ up to multiplication by a power of $t$.}

{\bf Proof}
Assume initially that ${\cal D}$ is the closure of a braid.

Write
$$M=\pmatrix{C_1& C_2&\ldots& C_n\cr}$$
in terms of its columns and let
$$M_{ij}=\pmatrix{C_1^0& C_2^0&\ldots& C_n^0\cr}$$
where each column has its $i$th component removed and $C_j^0$ does not
appear in the list. From $**$,
$$C_j^0=-C_1^0Q^{(1-j)}-\cdots-C_n^0Q^{n-j}$$
and $C_j^0$ does not
appear on the right hand side of the equation.

So by column operations which do not change the value of the
determinant we can change any column to $C_j^0$.
Now note that the value of the determinant is unchanged by
interchanging two columns. A similar argument works for the rows.

A general diagram is obtained from ${\cal D}$ by a sequence of
Reidemeister moves. A glance at the change of $M$ under the
Reidemeister moves shows that $\Delta_1$ is invariant up to
multiplication by a power of $t$ . \qed

Let us now return to $K_3$.
This is the closure of the braid 
$$\tau_2(\sigma_1\sigma_2\sigma_1)\tau_2
(\sigma_1\sigma_2\sigma_1)^{-1}$$
Suppose the representation of this as a $3\times 3$ matrix,
using the Budapest switch augmented by $t$, is $M$. Then the 
representation matrix of the module is $M-id$. The nine codimension
1 subdeterminants are
$$p=2t^{-2}+5+2t^{2},\hbox{ (4 times) } q=(2+2t^2)p(t^{-1}),
\hbox{ (twice) } q(t^{-1}), p^2$$
This not only shows that $K_3$ is non-trivial but that it cannot be classical
by 9.4.

\section{Epilogue}
Most calculations in this paper are done with Maple. At Andy Bartholomew's
website at 
http://www.layer8.co.uk/maths/
it is possible to download a C-program which calculates the invariants.

It is extremely unlikely that there are no non-trivial virtual knots for
which these methods fail to distinguish it from the trivial knot. For example
if the braid $\beta\in\overline{K}$ (see section 2), then the closure of 
$\beta$ possibly provides an infinite set of examples. However, to prove that
an infinite set exists would require different methods.

{\secfont References}
\refe

\cite{Alex} J. W. Alexander, Topological Invariants of Knots and Links,
Trans. American Math. Soc. 30(1928) pp 255-306

\cite{As} Helmer Aslaksen, Quaternionic Determinants, Math. Intel. Vol 18 no.
3 (1996)

\cite{BF} A. Bartholomew and Roger Fenn. Quaternionic Invariants of Virtual
Knots and Links, preprint. Preprint vailable from\nl
http://www.maths.sussex.ac.uk////Staff/RAF/Maths/Current/Andy/

\cite{BuF} S. Budden and Roger Fenn. The equation $$[b,(a-1)(a,b)]=0$$
and virtual knots and links, Fund Math 184 (2004) pp 19-29.

\cite{FJK} R. Fenn, M. Jordan, L. Kauffman,  Biquandles and 
Virtual Links, Topology and its Applications, 145 (2004) 157-175

\cite{CS} J. S. Carter, D. Jelsovsky, Seiichi Kamada,
Laurel Langford, Masahico Saito.
Quandle Cohomology and State-sum Invariants of Knotted Curves and Surfaces\nl
http://arxiv.org/abs/math.GT/9903135

\cite{De} P. Dehornoy, Non Commutative Versions of the Burau Representation,
C.R.Acad. Roy. Sci. Canada; 17-1; (1995) pp 53-58

\cite{Dr} V. Drinfeld, On some Unsolved Problems in Quantum Group Theory,
Quantum Groups, Lectures Notes in Maths. 1510, Springer 1-8 (1990)

\cite{ESS} P. Etingof, T. Schedler and A. Soloviev, Set-Theoretic Solutions to
the Quantum Yang-Baxter Equations. Duke Math Journal 100 no. 2 169-209 (1999)

\cite{EGS} P. Etingof, R. Guralnik and A. Soloviev, Indecomposable
set-theoretical solutions to the Quantum Yang-Baxter Equation on a set with
prime number of elements. J of Algebra 242 (2001) 709-719

\cite{FR} R. Fenn, C. Rourke. Racks and Links in Codimension Two. JKTR, No. 4,
343-406 (1992).

\cite{FRS1} R. Fenn, C. Rourke B. Sanderson. An Introduction to Species and the
Rack Space. M. E. Bozhuyu (ed.) Topics in Knot Theory, Kluwer Academic, pp
33-55.(1993)

\cite{FRS2} James Bundles (with C.Rourke and
B.Sanderson). {\sl Proceedings of the LMS} (3) 89 (2004) 217-240

\cite{J} D. Joyce A classifying invariant of knots, the knot quandle.
J. Pure Appl. Algebra 23, 37-65 (1982).

\cite{KK} N. Kamada and S. Kamada, Abstract Link Diagrams and Virtual
Knots,\nl JKTR 9, 93-106 (2000).

\cite{KK2} N. Kamada and S. Kamada, Braid Presentations of Virtual and
Welded Knots, preprint GT/0008092

\cite{Ku} G. Kuperberg, What is a Virtual Link? Algebraic and Geometric
Topology 587-591 (2003) 

\cite{K} L.Kauffman. Virtual Knot Theory, European J. Comb. Vol 20, 
663-690, (1999)

\cite{KS} T. Kishino and S. Satoh, A note on classical polynomials,
preprint (2001)

\cite{L} T. Y. Lam. The Algebraic Theory of Quadratic Forms, Benjamin (1973)

\bye

%% file: cprfonts.tex
\def\hexnumber#1{\ifcase#1 0\or 1\or 2\or 3\or 4\or 5\or 6\or 7\or 8\or
 9\or A\or B\or C\or D\or E\or F\fi}
%
%
\font\twelvemsa=msam10 scaled 1200   
\font\tenmsa=msam10                  
\font\ninemsa=msam9            \font\sevenmsa=msam7
\font\sixmsa=msam6             \font\fivemsa=msam5
%
%
\newfam\msafam                 \textfont\msafam=\tenmsa
\scriptfont\msafam=\sevenmsa   \scriptscriptfont\msafam=\fivemsa
\edef\hexa{\hexnumber\msafam}        
\def\msa{\fam\msafam\tenmsa}         
%
%
\font\twelvemsb=msbm10 scaled 1200   
\font\tenmsb=msbm10                  
\font\ninemsb=msbm9            \font\sevenmsb=msbm7
\font\sixmsb=msbm6             \font\fivemsb=msbm5
%
\newfam\msbfam                 \textfont\msbfam=\tenmsb       
\scriptfont\msbfam=\sevenmsb   \scriptscriptfont\msbfam=\fivemsb
\edef\hexb{\hexnumber\msbfam}        
\def\msb{\fam\msbfam\tenmsb}         
%
%
\font\twelveeufm=eufm10 scaled 1200  
\font\teneufm=eufm10                 
\font\nineeufm=eufm9           \font\seveneufm=eufm7
\font\sixeufm=eufm6            \font\fiveeufm=eufm5
%
\newfam\eufmfam                \textfont\eufmfam=\teneufm
\scriptfont\eufmfam=\seveneufm \scriptscriptfont\eufmfam=\fiveeufm
\edef\hexf{\hexnumber\eufmfam}      
\def\frak{\fam\eufmfam\teneufm}     
%
%
%
\font\twelverm=cmr10 scaled 1200    
\font\ninerm=cmr9                   
\font\sixrm=cmr6   
%
\font\twelvei=cmmi10 scaled 1200    
\font\ninei=cmmi9                   
\font\sixi=cmmi6  
%
\font\twelvesy=cmsy10 scaled 1200   
\font\ninesy=cmsy9                  
\font\sixsy=cmsy6  
%
\font\twelvebf=cmbx10 scaled 1200   
\font\ninebf=cmbx9                  
\font\sixbf=cmbx6  
%
%
\font\twelveit=cmti10 scaled 1200   
\font\nineit=cmti9                  
%
\font\twelvesl=cmsl10 scaled 1200   
\font\ninesl=cmsl9                  
%
\font\twelvett=cmtt10 scaled 1200   
\font\ninett=cmtt9                  
%
%
%
%
\def\small{%
%
%
\textfont0=\ninerm \scriptfont0=\sixrm \scriptscriptfont0=\fiverm
\def\rm{\fam0\ninerm}        
%
%
\textfont1=\ninei \scriptfont1=\sixi \scriptscriptfont1=\fivei
%
%
\textfont2=\ninesy \scriptfont2=\sixsy \scriptscriptfont2=\fivesy
%
%
\textfont3=\tenex \scriptfont3=\tenex \scriptscriptfont3=\tenex
%
%
\textfont\bffam=\ninebf \scriptfont\bffam=\sixbf
\scriptscriptfont\bffam=\fivebf \def\bf{\fam\bffam\ninebf}%
%
%
\textfont\itfam=\nineit \def\it{\fam\itfam\nineit}%
\textfont\slfam=\ninesl \def\sl{\fam\slfam\ninesl}%
\textfont\ttfam=\ninett \def\tt{\fam\ttfam\ninett}%
%
%
%
\textfont\msafam=\ninemsa \scriptfont\msafam=\sixmsa
\scriptscriptfont\msafam=\fivemsa \def\msa{\fam\msafam\ninemsa}%
%
%
\textfont\msbfam=\ninemsb \scriptfont\msbfam=\sixmsb
\scriptscriptfont\msbfam=\fivemsb \def\msb{\fam\msbfam\ninemsb}%
%
%
\textfont\eufmfam=\nineeufm  \scriptfont\eufmfam=\sixeufm
\scriptscriptfont\eufmfam=\fiveeufm \def\frak{\fam\eufmfam\nineeufm}%
%
%
%
\normalbaselineskip=11pt
\setbox\strutbox=\hbox{\vrule height8pt depth3pt width0pt}%
%
%
\normalbaselines\rm}    
%
%
%
%
\def\large{%
\textfont0=\twelverm \scriptfont0=\ninerm \scriptscriptfont0=\sevenrm
\def\rm{\fam0\twelverm}%
\textfont1=\twelvei \scriptfont1=\ninei \scriptscriptfont1=\seveni
\textfont2=\twelvesy \scriptfont2=\ninesy \scriptscriptfont2=\sevensy
\textfont3=\tenex \scriptfont3=\tenex \scriptscriptfont3=\tenex
\textfont\bffam=\twelvebf \scriptfont\bffam=\ninebf
\scriptscriptfont\bffam=\sevenbf \def\bf{\fam\bffam\twelvebf}%
\textfont\itfam=\twelveit \def\it{\fam\itfam\twelveit}%
\textfont\slfam=\twelvesl \def\sl{\fam\slfam\twelvesl}%
\textfont\ttfam=\twelvett \def\tt{\fam\ttfam\twelvett}%
\textfont\msafam=\twelvemsa \scriptfont\msafam=\ninemsa
\scriptscriptfont\msafam=\sevenmsa \def\msa{\fam\msafam\twelvemsa}         
\textfont\msbfam=\twelvemsb \scriptfont\msbfam=\ninemsb
\scriptscriptfont\msbfam=\sevenmsb \def\msb{\fam\msbfam\twelvemsb}         
\textfont\eufmfam=\twelveeufm  \scriptfont\eufmfam=\nineeufm
\scriptscriptfont\eufmfam=\seveneufm \def\frak{\fam\eufmfam\teneufm}
\normalbaselineskip=15pt
\setbox\strutbox=\hbox{\vrule height11pt depth4pt width0pt}%
\normalbaselines\rm}%
%
\def\Bbb{\msb}

%

%
\mathchardef\plussquare="0\hexa01
\mathchardef\nge="3\hexb0B
\mathchardef\maltesecross="0\hexa7A
\mathchardef\del="0\hexf01
%

%